\documentclass{gtmon_a}
\pdfoutput=1
\usepackage{amscd}
\input{diagxy}
\let\to\rightarrow
\let\barrsquare\square
\let\square\undefined

\proceedingstitle{Groups, homotopy and configuration spaces (Tokyo
  2005)}
\conferencestart{5 July 2005}
\conferenceend{11 July 2005}
\conferencename{Groups, homotopy and configuration spaces, 
                in honour of Fred Cohen's 60th birthday}
\conferencelocation{University of Tokyo, Japan}

\editor{Norio Iwase}
\givenname{Norio}
\surname{Iwase}

\editor{Toshitake Kohno}
\givenname{Toshitake}
\surname{Kohno}

\editor{Ran Levi}
\givenname{Ran}
\surname{Levi}

\editor{Dai Tamaki}
\givenname{Dai}
\surname{Tamaki}

\editor{Jie Wu}
\givenname{Jie}
\surname{Wu}

\title{String topology of Poincar\'e duality groups}

\author{Hossein Abbaspour}
\givenname{Hossein}
\surname{Abbaspour}
\address{Centre de Math\'ematiques\\
\'Ecole Polytechnique\\\newline
Palaiseau\\
France\vskip 3pt
Department of Mathematics\\
Stanford University\\\newline
Stanford CA 94305\\
USA\vskip 3pt
Department of Mathematics\\
Stanford University\\\newline
Stanford CA 94305\\
USA}

\email{abbaspou@mpim-bonn.mpg.de}
\urladdr{}

\author{Ralph Cohen}
\givenname{Ralph}
\surname{Cohen}
\email{ralph@math.stanford.edu}
\urladdr{}

\author{Kate Gruher}
\givenname{Kate}
\surname{Gruher}
\email{gruher@math.stanford.edu}
\urladdr{}

\dedicatory{Dedicated to Fred Cohen on the occasion of his 60th birthday}

\volumenumber{13}
\issuenumber{}
\publicationyear{2008}
\papernumber{1}
\startpage{1}
\endpage{10}

\doi{}
\MR{}
\Zbl{}

\arxivreference{math.AT/0511181}


\received{7 November 2005}
\revised{19 March 2007}
\accepted{22 March 2007}
\published{22 February 2008}
\publishedonline{22 February 2008}
\proposed{}
\seconded{}
\corresponding{}
\version{}

\makeatletter
\def\cnewtheorem#1[#2]#3{\newtheorem{#1}{#3}
\expandafter\let\csname c@#1\endcsname\c@theorem}

\newtheorem{theorem}{Theorem}
\cnewtheorem{lemma}[theorem]{Lemma}
\cnewtheorem{corollary}[theorem]{Corollary}
\cnewtheorem{proposition}[theorem]{Proposition}
\theoremstyle{remark}
\cnewtheorem{definition}[theorem]{Definition}
\cnewtheorem{remark}[theorem]{Remark}

\makeatother  

\let\xysavmatrix\xymatrix
\def\xymatrix{\disablesubscriptcorrection\xysavmatrix}
\AtBeginDocument{\let\bar\wbar\let\tilde\widetilde\let\hat\what}

\makeop{ev}
\makeop{Ad}

\def\lbracket{[}
\def\rbracket{]}
\def\semicolon{;}

\newcommand{\la}{\to}
\newcommand{\hk}{\hookrightarrow}

\newcommand{\br}{{\mathbb R}}
\newcommand{\bz}{{\mathbb Z}}

\newcommand{\bh}{{\mathbb H}}

\newcommand{\cll}{\mathcal{L}}

\newcommand{\xr}{\xrightarrow}
\newcommand{\bfl}{\begin{flushleft}}
\newcommand{\efl}{\end{flushleft}}

\begin{document}

\begin{htmlabstract}
Let G be a Poincar&eacute; duality group of dimension n.  For a given
element g &isin; G, let C<sub>g</sub> denote its centralizer
subgroup. Let L<sub>G</sub> be the graded abelian group defined by
(L<sub>G</sub>)<sub>p</sub> =
&oplus;<sub>[g]</sub>H<sub>p+n</sub>(C<sub>g</sub>)
where the sum is taken over conjugacy classes of elements in G.  In
this paper we construct a multiplication on L<sub>G</sub> directly in
terms of intersection products on the centralizers.  This
multiplication makes L<sub>G</sub> a graded, associative, commutative
algebra.  When G is the fundamental group of an aspherical, closed
oriented n&ndash;manifold M, then (L<sub>G</sub>)<sub>*</sub> =
H<sub>*+n</sub>(LM), where LM is the free loop space of M.  We show
that the product on L<sub>G</sub> corresponds to the string topology
loop product on H<sub>*</sub>(LM) defined by Chas and Sullivan.
\end{htmlabstract}

\begin{abstract}
Let $G$ be a Poincar\'e duality group of dimension $n$.  For a given
element $g \in G$, let $C_g$ denote its centralizer subgroup. Let
$\mathcal{L}_G$ be the graded abelian group defined by
$$(\mathcal{L}_G)_p = {\textstyle\bigoplus_{[g]}}H_{p+n}(C_g)$$
where the sum is taken over conjugacy classes of elements in $G$.  In this
paper we  construct a multiplication on $\mathcal{L}_G$ directly in terms
of intersection products on the centralizers.  This multiplication makes
$\mathcal{L}_G$ a graded, associative, commutative algebra.  When $G$ is
the fundamental group of an aspherical, closed oriented $n$--manifold $M$,
then $(\mathcal{L}_G)_* = H_{*+n}(LM)$, where $LM$ is the free loop space
of $M$.  We show that the product on $\mathcal{L}_G$ corresponds to the
string topology loop product on $H_*(LM)$ defined by Chas and Sullivan.
\end{abstract}

\maketitle

\section*{Introduction}\addcontentsline{toc}{section}{Introduction}
In \cite{chassullivan} Chas and Sullivan defined an associative,
commutative algebra structure on $H_*(LM)$, the homology of the free
loop space of a closed, oriented $n$--manifold $M$. The Chas--Sullivan
string topology product has total degree $-n$, $H_p(LM) \otimes H_q(LM)
\to H_{p+q-n}(LM)$.

On the other hand it is well known that if $G$ is any topological group,
and $BG$ is its classifying space, then its loop space is naturally homotopy equivalent
to the homotopy orbit space of the adjoint action of the group on itself:
$$L(BG) \simeq EG \times_G G$$
where $EG$ is a contractible space with a free $G$ action, and $G$ acts on
itself by conjugation.  The space $EG \times_G G$ is the orbit space of
the corresponding diagonal $G$--action on $EG \times G$. 
When $G$ is discrete, we therefore have an isomorphism,
$$H_*(L(BG)) \cong H_*(G; \Ad)$$
where $\Ad$ is the $G$--module given by $\bz [G]$ additively, and the action is given by conjugation. It is an elementary fact that this homology group splits as
$$H_*(G; \Ad) = \oplus_{[g]} H_*(C_g)$$
where the direct sum is taken over conjugacy classes of elements in $G$, and $C_g$ is the centralizer of  $g \in G$.  

Combining this with Chas and Sullivan's string topology construction, when $G$ is the fundamental group of a closed, oriented, aspherical manifold $M$, then
$$(\cll_G)_* := \oplus_{[g]} H_{*+n}(C_g)$$
has the structure of a graded commutative algebra.

It is the goal of this note to describe this product structure directly
in terms of homology of groups. Our construction will not involve any
differential topology, and so will work for any \emph{Poincar\'e duality
group}: that is a group $G$ together with a fundamental class $z \in
H_n(G; \bz)$ such that capping with this class is an isomorphism
\begin{equation}\label{pd}
\begin{CD}
\cap z \co H^p(G, A) @>\cong >>  H_{n-p}(G, A)
\end{CD}
\end{equation}
where $A$  is an arbitrary $\bz [G]$--module.  Such a structure is called
an oriented Poincar\'e duality group in Brown \cite{brown}.

Our construction will be defined in terms of the \emph{intersection
product} in the homology of subgroups of a Poincar\'e duality group.
Given two such subgroups, $K, H < G$, this product is a pairing,
$$
H_p(K; \bz) \otimes H_q(H; \bz) \la \bigoplus_{[g]} H_{p+q-n}(K \cap gHg^{-1}).
$$
The target of this map is the sum of the homologies of the intersections of K with conjugates of $H$,    taken over all double cosets $[g] = KgH \subset G$.  
This construction is defined directly in terms of the Poincar\'e duality isomorphism \eqref{pd}.
In the case when $K$ and $H$ are centralizers, we will be able to assemble the intersection maps into a product on $\cll_G$.  After making this construction we will show that this operation corresponds to the  
Chas--Sullivan string topology loop  product in the case when $G$ is the fundamental group of an aspherical manifold (see \fullref{main} below).

\section{Intersection pairings in Poincar\'e duality groups and the string product}

Throughout this section $G$ will be a fixed Poincar\'e duality group of
dimension $n$.  Let $z \in H_n(G)$ be a fundamental class.  We let $D\co
H^q(G; M) \smash{\xr{\cong}} H_{n-q}(G;M)$ be the Poincar\'e duality isomorphism
given by capping on the left with $z$.

We begin by using the Poincar\'e duality isomorphism to define an
\emph{intersection pairing},
$$\cap \co H_p(K; \bz) \otimes H_q(H; \bz) \la \bigoplus_{[g]\in I}
H_{p+q-n}\bigl(K \cap gHg^{-1}\bigr),$$
where $K$ and $H$ are any two subgroups of $G$, and $I$ is an indexing
set for the double cosets $K\backslash G /H$.

First  observe that for any  subgroup $K$ of $G$ we have
$$H_*(K,\mathbb{Z})=H_*(G,\mathbb{Z}[G/K])$$
where $\mathbb{Z}[G/K]$ is the free $G$--module generated by the left
cosets, equipped with the standard left action of $G$.

We next need the following double coset formula.

\begin{lemma} \label{doublecoset}
Suppose $I$ is a set of representatives for the double cosets $K\backslash
G /H$. Then there is an isomorphism of $\mathbb{Z}[G]$--modules,
\begin{equation}
\mathbb{Z}[G/K]\otimes \mathbb{Z}[G/H]\cong \bigoplus_{g\in I}
\mathbb{Z}\bigl[G/K\cap gH g^{-1}\bigr],
\end{equation}
where $G$ acts diagonally on $\mathbb{Z}[G/K]\otimes \mathbb{Z}[G/K]$
by the standard left translation.
\end{lemma}

\begin{proof}
Given $g_1K\otimes g_2H$,  double coset $K(g_1^{-1}g_2)H$ is well
defined and does not depend on the choice of the representatives
$g_1$ and $g_2$ for the cosets $g_1K$ and $g_2H$, so let $g\in I$ be
representative for the double $K(g_1^{-1}g_2)H$. If we write $g_1K\otimes
g_2H=a(K\otimes gH)$ for some  $a\in G$, the isomorphism is defined
by sending $g_1K\otimes g_2H$ to  $a(K\cap gH g^{-1})$. This map is
well-defined as $a(K\otimes gH)=b(K\otimes gH)$ implies  that $b^{-1}a\in
K\cap gH g^{-1}$.  It is obvious that the map described above is a map
of $G$--modules.
\end{proof}

With this we can define the intersection product.

\begin{definition}\label{intersect}
The intersection pairing
$$\cap \co  H_i(K; \bz) \otimes H_j(H; \bz) \la \bigoplus_{[g]} H_{i+j-n}(K
  \cap gHg^{-1})$$
is given by the composition
\begin{multline*}
\cap \co H_i(K) \otimes H_j(H) \xr{~\cong~} H_i(G,\, \mathbb{Z}[G/K])\otimes
  H_j(G, \, \mathbb{Z}[G/H]) \xr{D^{-1} \otimes D^{-1}} \\
H^{n-i}(G,\mathbb{Z}[G/K])\otimes H^{n-j}(G,\mathbb{Z}[G/H])
  \xr{\text{cup prod.}} H^{2n-i-j}(G,\mathbb{Z}[G/K]\otimes
  \mathbb{Z}[G/H])\\
\xr{~D~} H_{i+j-n}(G,\mathbb{Z}[G/K]\otimes \mathbb{Z}[G/H]) \xr{~\psi_*~}
  \bigoplus_{[g]\in I} H_{i+j-n}(G; \bz[G/K\cap gH g^{-1}]) \\
\xr{~\cong~} \bigoplus_{[g]\in I} H_{i+j-n}(K \cap gHg^{-1})
\end{multline*}
where the map $\psi$ is the isomorphism given by \fullref{doublecoset}.
\end{definition}

\begin{remark}
When $K = H = G$, the intersection pairing
$$\cap \co H_i(G) \otimes H_j(G) \to  H_{i+j-n}(G)$$
is simply the homomorphism that is Poincar\'e dual to the cup product in cohomology.
In the case when $G$ is the fundamental group of a closed, oriented,
aspherical $n$--manifold $M$, this is the usual definition (up to sign) of the intersection product on $H_*(M)$.
\end{remark}

The case of interest for us is when the subgroups are the centralizers
of the elements of $G$. Suppose $C_{\alpha}$ and $C_{\beta}$ are the
centralizers of $\alpha$ and  $\beta$ in $G$.  Then we have
\begin{equation}
\cap\co H_i(C_{\alpha})\otimes H_j(C_{\beta}) \rightarrow
\bigoplus_{g\in I}H_{n-i-j}(G,\mathbb{Z}[G/C_{\alpha}\cap
gC_{\beta} g^{-1}]).
\end{equation}
Now the conjugate of a centralizer subgroup is the centralizer of the conjugate,
$$gC_{\beta}g^{-1} = C_{g \beta g^{-1}}.$$
Furthermore we have inclusions of subgroups, 
$$j \co C_{\alpha}\cap gC_{\beta}g^{-1}=C_{\alpha}\cap C_{g
  \beta g^{-1}} \leq C_{\alpha g\beta g^{-1}}.$$
This will allow us to make the following definition.

\begin{definition}\label{loopprod}
Define the string product,
\begin{align*}
\mu \co (\cll_G)_p  \otimes  (\cll_G)_q &\longrightarrow (\cll_G)_{p+q} \\
\Bigl( \bigoplus_{\alpha} H_{p+n}(C_\alpha)\Bigr) \otimes
\Bigl(\bigoplus_{\beta} H_{q+n}(C_\beta) \Bigr) &\longmapsto
\Bigl(\bigoplus_{\gamma} H_{p+q+n}(C_\gamma)\Bigr)
\end{align*}
to be the composition
\begin{multline*}
\mu = j_* \circ \cap \co H_{p+n}(C_\alpha) \otimes H_{q+n}(C_\beta)
\longrightarrow \\
\bigoplus_{g\in I}H_{p+q+n}(C_\alpha \cap
gC_{\beta}g^{-1})  \longrightarrow \bigoplus_{g\in I}H_{p+q+n}(C_{\alpha g\beta g^{-1}}).
\end{multline*}
\end{definition}

\begin{theorem}
The product $\mu$ makes $\cll_G$ a graded,  associative, commutative
algebra  with unit.
\end{theorem}

\begin{proof}
The associativity is a direct check of the definitions.  The unit
in the algebra is given by the element $z \in H_n(G) \subset
(\cll_G)_0$. Again, it is an immediate check of definitions to verify
that this class acts as the identity.

To see commutativity, consider the map
$$\bigoplus_{\alpha, \beta} \bigoplus_{g\in I}H_{*+n}( C_{\alpha}\cap
  C_{g{\beta} g^{-1}})  \xr{j_*} \bigoplus_{\alpha, \beta} \bigoplus_{g \in
  I} H_{*+n}(C_{\alpha g\beta  g^{-1}}) \to  (\cll_G)_* = H_{*+n}(G, \Ad).$$
Notice that $ \bigoplus_{\alpha, \beta} \bigoplus_{g\in I}H_{*+n}(
C_{\alpha}\cap C_{g{\beta} g^{-1}})$ has an involution on it, by
interchanging the roles of $\alpha$ and $\beta$, and mapping $g$
to $g^{-1}$.   We call that involution $\tau$.  Notice that for $x \in
H_{p+n}(C_{\alpha})$, and $y \in H_{q+n}(C_\beta)$, we have that
$$\mu (y,x) = j_* \circ \cap (y, x) = (1)^{pq} j_* \circ \tau \circ \cap
(x,y).$$
Now notice  the obvious conjugation relation,
$$g^{-1}\alpha^{-1}(C_{\alpha g\beta g^{-1}})\alpha g = C_{\beta g^{-1}
  \alpha g},$$
Using the fact that if $h_1$ is conjugate to $h_2$ then the images of the
homology of the centralizers, $H_*(C_{h_1}) \hk H_*(G; \Ad)$ and
$H_*(C_{h_2}) \hk H_*(G; \Ad)$ are equal, we see that $j_* \circ \tau = j_*$, and so
$\mu (x,y) = (-1)^{pq} \mu (y,x)$. 
\end{proof}

Now as mentioned above, the free loop space, $LBG$ is homotopy equivalent
to the homotopy orbit space of the adjoint action, $EG \times_{\Ad} G$.  This gives a natural isomorphism of homologies,
$$H_*(LBG) \cong H_*(G; \Ad).$$
Now assume that $G = \pi_1(M^n)$, where $M^n$ is an aspherical, closed,
oriented $n$--manifold.  In particular this says that $M^n \simeq BG$, and so there is a homotopy equivalence,
$$LM \simeq EG \times_{\Ad} G$$
and so an isomorphism,
$$H_*(LM) \cong H_*(G; \Ad).$$
Using transversal intersection theory, in \cite{chassullivan} Chas and Sullivan described a string topology loop product on the regraded homology groups,
$\bh_*(LM^n) = H_{*+n}(LM)$,
$$\circ \co \bh_*(LM^n) \otimes  \bh_*(LM^n) \la  \bh_*(LM^n).$$
 Notice with respect to this regrading, we have an isomorphism,
$\Phi\co \bh_*(LM^n) \cong ( \cll _G)_*$.  The following says that
our algebraically defined product on $(\cll_G)_*$ corresponds to the
Chas--Sullivan string topology product on $\bh_*(LM)$.

\begin{theorem}\label{main}
The following diagram commutes:
$$\bfig
  \barrsquare<1000,500>[\bh_*(LM^n) \otimes  \bh_*(LM^n)`
    \bh_*(LM^n)`  (\cll_G)_*  \otimes (\cll_G)_*`  (\cll_G)_*;
    \circ`\Phi`\Phi`\mu]
  \morphism(0,0)|r|//<0,500>[\phantom{\bh_*(LM^n) \otimes  \bh_*(LM^n)}`
    \phantom{(\cll_G)_*  \otimes (\cll_G)_*};\cong]
  \morphism(1000,0)|l|//<0,500>[\phantom{\bh_*(LM^n)}`
    \phantom{(\cll_G)_*};\cong]
  \efig$$
\end{theorem}

\begin{proof}
Recall that the string topology product was described by Cohen and Jones
\cite{cohenjones} in the following way.  Consider the diagram
$$\begin{CD}
  LM @<\gamma << LM \times_{M} LM @>\iota >> LM \times LM
  \end{CD}$$
where $LM \times_{M} LM$ is the subspace of $LM \times LM$ 
consisting of those pairs of loops $(\alpha, \beta)$, with $\alpha (0) =
\beta (0)$.  $\iota \co  LM \times_{M} LM \hk LM \times LM $ is the natural
inclusion.  $\gamma \co  LM \times_{M} LM \to LM$ is the map obtained by
thinking of $ LM \times_{M} LM$ as the maps from the figure 8 to M, $ LM
\times_{M} LM \cong Map(8, M)$, and then mapping this space to $LM$ by
sending a map from the figure 8 to the loop obtained by starting at the
intersection point of the 8, first traversing the upper circle, and then
traversing the lower circle of the 8.  See Chas--Sullivan
\cite{chassullivan} and Cohen--Jones \cite{cohenjones} for details.   Notice that $ LM \times_{M} LM$    is the fiber product for the square,
\begin{equation}\label{pullback}
\bfig
\barrsquare<1000,400>[LM\times_{M}LM` LM \times LM` M` M \times M;
  \iota``\ev \times \ev`\Delta]
\efig
\end{equation}
Here $\ev \co LM \to M$ evaluates a loop at the starting point, $0 \in
\br/\bz = S^1$.  In \cite{cohenjones} a Pontrjagin--Thom map $\tau_\iota
\co LM \times  LM \to  (LM \times_{M} LM)^{TM}$ was constructed, where the
target is the Thom space of the tangent bundle $TM \to M$ pulled back
over $ LM \times_{M} LM$.  This map was defined by finding a tubular
neighborhood of the embedding $\iota$, and collapsing everything in $LM
\times LM$ outside the tubular neighborhood to a point.  When one applies
homology and the Thom isomorphism to the map $\tau_\iota$ one obtains an
``umkehr map",
\begin{equation}\label{umkehr}
\iota_! \co H_*(LM \times LM) \to H_{*-n}( LM \times_{M} LM)
\end{equation}
and the basic string topology product is the composition,
$$\circ = \gamma_* \circ \iota_! \co  H_*(LM \times LM) \to H_{*-n}(
LM \times_{M} LM) \to H_{*-n}(LM).$$
See \cite{cohenjones} for details. 

Now in the case of interest, where $M$ is an aspherical manifold
with fundamental group $G$,  the tubular neighborhood of $\iota\co
LM \times_{M} LM \hk LM \times LM$, and the resulting Thom--Pontrjagin
construction has the following description.

Let $\tilde M \to M$ be the universal cover.  Since $M$ is aspherical,
$\tilde M$ is a contractible space with a free $G$--action.   We can
therefore write the homotopy equivalence,
$$LM \cong \tilde M \times_{\Ad} G,$$
and up to homotopy, the evaluation map $\ev \co LM \to M$ is given by the
projection map $p\co \tilde M \times_{\Ad} G \to M.$  Now as seen before,
$\tilde M \times_{\Ad} G$ decomposes as a disjoint union of spaces,
$$\tilde M \times_{\Ad} G = \coprod_{[\alpha]} \tilde M \times_{G}
(G/C_\alpha)$$
where the disjoint union is taken over conjugacy classes of elements
$\alpha \in G$. Thus the pullback square \eqref{pullback} is the disjoint
union, over pairs of conjugacy classes, $\alpha, \beta$, of the pullback
squares,
\begin{equation}\label{grouppull}
\bfig
\barrsquare<1500,500>[\tilde M{\times_G}(G/C_\alpha{\times}G/C_{\beta})`
  (\tilde M{\times}\tilde M){\times_{G\times G}}(G/C_\alpha{\times}G/C_{\beta})`
  M` M{\times}M;
  \tilde \Delta{\times}1` p` p{\times}p` \Delta]
\efig
\end{equation}
A tubular neighborhood of the diagonal map, $\Delta \co M \to M\times M$
lifts to a $G$--equivariant tubular neighborhood  of the diagonal, $\tilde
\Delta\co \tilde M \to \tilde M \times \tilde M,$  which has normal bundle
isomorphic to the tangent bundle, $T\tilde M \to \tilde M$.  Therefore the
Pontrjagin--Thom collapse map in this setting, is $G$--equivariant,
$$\tau \co \tilde M \times \tilde M \to  (\tilde M)^{T\tilde M}.$$
This says that the Pontrjagin--Thom map for the embedding $\tilde \Delta
\times 1$ from diagram \eqref{grouppull}, a disjoint union of which
model the Pontrjagin--Thom map of the embedding $\iota$ from diagram
\eqref{pullback}, is given by
$$\tau \times 1 \co (\tilde M \times \tilde M) \times_{G\times G}
(G/C_\alpha \times G/C_{\beta}) \la  (\tilde M)^{T\tilde M}_+
\wedge_{G\times G} (G/C_\alpha \times G/C_{\beta}).$$
Applying homology and the Thom isomorphism, one gets  umkehr maps,
$$(\iota_{\alpha, \beta})_! \co H_*((\tilde M \times \tilde M)
\times_{G\times G} (G/C_\alpha \times G/C_{\beta})) \la H_{*-n}(\tilde
M \times_{G} (G/C_\alpha \times G/C_{\beta})),$$
A direct sum of these umkehr maps, (over conjugacy classes $\alpha$
and $\beta$) gives a model for the umkehr map $\iota_!$ given in
\eqref{umkehr} above:
\begin{equation}\label{commute}
\begin{footnotesize}
  \bfig
  \barrsquare<1800,500>[\bigoplus_{\alpha,\beta}H_*((\tilde M{\times}\tilde M)
  {\times}_{G{\times}G}(G/C_\alpha {\times} G/C_{\beta}))`
  \bigoplus_{\alpha, \beta}H_{*-n}(\tilde M {\times}_{G}
  (G/C_\alpha{\times} G/C_{\beta}))`
  H_*(LM {\times} LM)`  H_{*-n}(LM{\times}_M LM);
  \bigoplus(\iota_{\alpha,\beta})_{!}`\cong`\cong`\iota_!]
  \efig
\end{footnotesize}
\end{equation}
In terms of group homology, the top line of this diagram is a sum of terms of the form,
$$(\iota_{\alpha, \beta})_{!}\co  H_*(G\times G; \,  \bz[G/C_\alpha \times
  G/C_\beta])  \la H_{*-n}(G; \,  \bz[G/C_\alpha \times G/C_\beta]).$$
Now the standard relationship between umkehr maps for embeddings of closed manifolds and Poincar\'e duality says that $(\iota_{\alpha, \beta})_!$ is Poincar\'e dual to the induced homomorphism in cohomology,
$$(\tilde \Delta \times 1)^* \co H^*((\tilde M \times \tilde M)
\times_{G\times G} (G/C_\alpha \times G/C_{\beta})   \la H^*(\tilde M
\times_G (G/C_\alpha \times G/C_{\beta}))$$
or, in terms of group cohomology,
$$(\tilde \Delta \times 1)^* \co H^*((G\times G; \, \bz[G/C_\alpha \times
G/C_\beta])  \la H^*(G; \,  \bz[G/C_\alpha \times G/C_\beta]).$$
But the diagonal map in group cohomology induces the cup product,
$$\cup \co H^*(G; \, \bz [G/C_\alpha]) \otimes H^*(G; \, \bz[G/C_\beta])
  \to H^*(G; \,  \bz[G/C_\alpha \times G/C_\beta]).$$
Putting these observations together with the definition of the intersection pairing in the homologies of the centralizers (\fullref{loopprod}), we see that the following diagram commutes: 
$$\bfig
  \barrsquare/->`->``->/<1800,800>[
    \bigoplus_{\alpha,\beta} H_*(C_\alpha) \otimes H_*(C_\beta)`
    \bigoplus_{\substack{\alpha,\beta \\ \lbracket g\rbracket \in I}}
    H_{*-n}(C_\alpha \cap gC_\beta g^{-1})`
    H_*(LM) \otimes H_*(LM)`  H_{*-n}(LM \times_M LM);
    \cap` \cong` ` \iota_!]
  \morphism(1800,800)<0,-400>[\phantom{H_{*-n}(C_\alpha \cap gC_\beta g^{-1})}`
    \bigoplus_{\alpha, \beta} H_{*-n}(G\semicolon\,\bz\lbracket
    G/C_{\alpha}{\times}G/C_\beta\rbracket);\cong]
  \morphism(1800,400)<0,-400>[\phantom{\bigoplus_{\alpha, \beta}
    H_{*-n}(G\semicolon\,\bz\lbracket
    G/C_{\alpha}{\times}G/C_\beta\rbracket)}`
    \phantom{H_*(LM) \otimes H_*(LM)`  H_{*-n}(LM \times_M LM)};\cong]
  \efig$$
On the other hand, it is clear that the following diagram also commutes,
$$\bfig
  \barrsquare/->`->``->/<1800,800>[
    \bigoplus_{\alpha, \beta}\bigoplus_{\lbracket g \rbracket \in I}
      H_{*-n}(C_\alpha \cap gC_\beta g^{-1})`
    \bigoplus_{\alpha, \beta}\bigoplus_{\lbracket g \rbracket \in I}
      H_{*-n}(C_{\alpha g\beta g^{-1}})`
    H_{*-n}(LM \times_M LM)` H_{*-n}(LM);
    j`\cong``\gamma_*]
  \morphism(1800,800)/^{ (}->/<0,-400>[\phantom{\bigoplus_{\alpha,
    \beta}\bigoplus_{\lbracket g \rbracket \in I} H_{*-n}(C_{\alpha g\beta
    g^{-1}})}` (\cll_G)_*;]
  \morphism(1800,400)<0,-400>[\phantom{(\cll_G)_*}`
    \phantom{H_{*-n}(LM \times_M LM)` H_{*-n}(LM)};\cong]
  \efig$$
The composition of the bottom rows of these two diagrams is the string topology loop product.  The composition of the top rows of these two diagrams is the map $\mu$.  Thus the theorem is proven.
\end{proof}

\subsection*{Acknowledgements}\addcontentsline{toc}{subsection}{Acknowledgements}
The first author was  supported by a ChateauBriand Postdoctoral
Fellowship.  The second author was partially supported by a research
grant from the NSF, and the third author was supported by a National
Defense Science and Engineering Graduate Fellowship.

\bibliographystyle{gtart}
\bibliography{link}

\end{document}